\documentclass[a4paper,11pt]{amsart}

\usepackage{amsfonts,latexsym,rawfonts,amsmath,amssymb,amsthm, mathrsfs, lscape}
\usepackage[all]{xy}
\usepackage[english]{babel}
\usepackage[utf8]{inputenc}

\usepackage{graphicx}

\usepackage{xcolor}

\usepackage{tikz-cd}

\usepackage{array, tabularx}

\usepackage{setspace}
\usepackage{textcomp}

\usepackage[hypertexnames=false,
backref=page,
    pdftex,
    pdfpagemode=UseNone,
    breaklinks=true,
    extension=pdf,
    colorlinks=true,
    linkcolor=blue,
    citecolor=blue,
    urlcolor=blue,
]{hyperref}

\usepackage[justification=centering]{caption}

\usepackage[top=1.5in, bottom=.9in, left=0.9in, right=0.9in]{geometry}

\newcommand{\referenza}{}

\newtheorem{thm}{Theorem}[section]
\newtheorem*{thm*}{Theorem \referenza}
\newtheorem{cor}[thm]{Corollary}
\newtheorem*{cor*}{Corollary \referenza}

\newtheorem*{lem*}{Lemma \referenza}

\newtheorem{prop}[thm]{Proposition}
\newtheorem*{prop*}{Proposition \referenza}

\newtheorem{conj}[thm]{Conjecture}
\newtheorem*{conj*}{Conjecture \referenza}

\theoremstyle{remark}
\newtheorem{rmk}[thm]{Remark}
\newtheorem*{rmk*}{Remark}

\newtheorem{exa}[thm]{Example}

\theoremstyle{definition}
\newtheorem{defi}[thm]{Definition}

\numberwithin{equation}{section}

\newcommand*{\DashedArrow}[1][]{\mathbin{\tikz [baseline=-0.25ex,-latex, dashed,#1] \draw [#1] (0pt,0.5ex) -- (1.3em,0.5ex);}}
\newcommand{\rat}{\DashedArrow[->,densely dashed]}

\def \O {\mathcal O}

\def \H {\mathbb H}

\def \S {\mathbb S}
\def \C {\mathbb C}
\def \Z {\mathbb Z}
\def \P {\mathbb P}

\def \p {\partial}

\renewcommand{\p}{\partial}
\renewcommand{\bar}{\overline}

\allowdisplaybreaks[1]

\title[]{A survey on rational curves on complex surfaces}

\author{Giuseppe Barbaro}
\address[Giuseppe Barbaro]{Dipartimento di Matematica ``Guido Castelnuovo", Sapienza Università di Roma, Piazzale Aldo Moro~5, 00185 Roma, Italy}
\email{g.barbaro@uniroma1.it}

\author{Filippo Fagioli}
\address[Filippo Fagioli]{Dipartimento di Matematica e Informatica ``Ulisse Dini", Università degli Studi di Firenze, Viale Morgagni~67/a, 50134 Firenze, Italy}
\email{filippo.fagioli@unifi.it}

\author{Ángel David Ríos Ortiz}
\address[Ángel David Ríos Ortiz]{Max-Planck-Institut f\"{u}r Mathematik in den Naturwissenschaften, Inselstrasse~22, 04103 Leipzig, Germany}
\email{arios@mis.mpg.de}

\thanks{The authors are supported by GNSAGA of INdAM. The first and second authors are also supported by project PRIN2017 ``Real and Complex Manifolds: Topology, Geometry and holomorphic dynamics'' (code 2017JZ2SW5)} 

\begin{document}

\begin{abstract}
    In this survey, we discuss the problem of the existence of rational curves on complex surfaces, both in the K\"ahler and non-K\"ahler setup. We systematically go through the Enriques--Kodaira classification of complex surfaces to highlight the different approaches applied to the study of rational curves in each class. We also provide several examples and point out some open problems.
\end{abstract}
	
\maketitle
	
\section{Introduction}
    The study of rational curves\footnote{In the following, by a rational curve in a compact complex manifold $X$ we mean a non-constant holomorphic map $\P^1\to X$, and we will usually identify the map with its image in $X$.} on complex manifolds has proved to be an important topic in Complex Geometry from the very beginning. For complex surfaces, the presence of rational curves affects directly the geometry of the surface. 
	This is traced back to the Italian school of Algebraic Geometry in the works of Castelnuovo and  Enriques in the classification of algebraic surfaces and to the work of Kodaira on the existence of minimal models for complex surfaces.
	Thanks to the {\em Castelnuovo Contraction Theorem} (and its generalization to all complex surfaces by Kodaira) the rational curves on surfaces can be studied on their {\em minimal models}. 
	Indeed, this result ensures that if a non-singular surface $S$ has a curve $C$ of self-intersection $-1$ (called \emph{$(-1)$-curve}), then there is a holomorphic map to another non-singular surface which contracts $ C $ to a point, and is a biholomorphism outside $ C $; the iteration of this process of contractions leads to a minimal model of $S$.
	Accordingly, we analyze the existence of rational curves on minimal surfaces (i.e. free from $(-1)$-curves), following the Enriques--Kodaira classification as in Table~\ref{tab: Enr-Kod classification}.
	
	\begin{table}[h!] 
		\begin{tabular} {|p{7cm}|c|c|c|c|}
			\hline 
			&&&&\\
			Class of $S$                  &$\operatorname{Kod}(S)$ &$a(S)$&$\chi(S)$&$b_1(S)$ \\
			&&&&\\  \hline  &&&&\\
			1) Rational surfaces          &               & 2    &  3,4    &  0    \\
			2) Ruled surfaces of genus $g\geq 1$&$-\infty$& 2    &$4(1-g)$ &  2g \\
			3) Class VII$_0^0$ surfaces     &               & 0,1  &  0      &  1    \\
			4) Class VII$_0^+$ surfaces   &               & 0    &  $>0$   &  1    \\
			&&&&\\    \hline    &&&&\\
			5) K3 surfaces                &         & 0,1,2& 24      & 0     \\
			6) Enriques surfaces          &         & 2    & 12      & 0     \\
			7) Complex tori                       &    0     & 0,1,2& 0       & 4     \\
			8) Hyperelliptic surfaces     &        & 2    & 0       & 2     \\
			9) Kodaira surfaces           &         & 1    & 0       & 1,3   \\
			&&&& \\ \hline &&&& \\
			10) Properly elliptic surfaces & 1       & 1,2    &$\geq 0$ & even  \\
			&         & 1      & 0       & odd   \\
			&&&& \\ \hline &&&& \\
			11) Surfaces of general type  & 2       & 2    & $>0$    & even  \\
			&&&& \\\hline
		\end{tabular}
		\vspace{1ex}
		\caption{Enriques--Kodaira classification}
		\label{tab: Enr-Kod classification}
	\end{table}
	
	The existence of rational curves is crucial for the classification of class VII surfaces. 
	These are non-K{\"a}hler compact complex surfaces with Kodaira dimension $-\infty$, and the minimal ones are denoted in the literature as class VII$_0$ surfaces. 
	The class VII$_{0}$ surfaces with second Betti number $b_2=0$ (denoted by class VII$_0^0$) have been classified \cite{Bog76, Bog82, LYZ90, LYZ94, Tel94}, and are either Hopf or Inoue--Bombieri surfaces. 
    On the other hand, those with $b_2>0$ (which are said class VII$_0^+$ or minimal class VII$^+$) are not classified in general yet. 
    Indeed, only those with $b_2=1$ were classified.
    This is due to Nakamura, who did it in \cite{Nak84} under the additional assumption that the surface has a curve, and to Teleman, who proved this extra condition in \cite{Tel05}.
    The dedicated effort of many authors, culminating in the theorem of Dloussky--Oeljeklaus--Toma \cite{DOT}, has reduced the problem of classifying class VII$^+$ surfaces to finding $b_2$ rational curves in a given minimal surface of class VII$^+$; this problem is known in the literature as {\em Global Spherical Shell Conjecture} (see \cite{Nak84b} and the references therein).
	
	\medskip
	
	The interest to write this survey arose as an attempt of the authors to put together, in a  systematic form, the state of the art on the existence of rational curves on complex surfaces. We summarize the collected results in the Table~\ref{table: summary on curves over surfaces}, where we specify whether there are finitely many, infinitely many or no rational curves.
	
	\begin{table}[h!]
		\begin{tabular} {|c|p{5cm}|c|}
			\hline 
			&&\\
			$\operatorname{Kod}(S)$ &Class of $S$& Rational curves on $ S $ \\
			&&\\  \hline  &&\\
			&1) Rational surfaces&covered by\\
			&2) Ruled surfaces of $g\geq 1$&covered by\\
			$-\infty$ &3) Class VII$_0$ surfaces& (Table~\ref{table: class VII}) \\
			&&\\    \hline    &&\\
			&4) K3 surfaces&no (Remark~\ref{rem:K3withnocurves})\\
			& &finitely many (Example~\ref{exa:kummer})\\
			& &infinitely many (Theorem~\ref{thm:infinityrationalcurvesK3})\\
			0&5) Enriques surfaces&infinitely many (Theorem~\ref{thm:Enriquesinfinitecurves})\\
			&6) Complex tori&no (Proposition~\ref{prop:norationalcurvestorus})\\
			&7) Hyperelliptic surfaces&no (Corollary~\ref{cor:norationalhyperelliptic})\\
			&8) Kodaira surfaces&no (Propositions~\ref{prop: primary kodaira surfaces case} and~\ref{prop: secondary kodaira surfaces case})\\
			&& \\ \hline && \\
			1&9) Properly elliptic surfaces&  \\
			&\qquad non-K{\" a}hler&no (Theorem~\ref{thm:nonkahlerproperlyelliptic})\\
			&\qquad K{\" a}hler&no (Example~\ref{ex: elliptic algebraic no curves})\\
			& &finitely many (Theorem~\ref{thm: ulmer theorem})\\
			& &infinitely many (Remark~\ref{rmk: infinite curves on elliptic})\\
			&& \\ \hline && \\
			2&10) Surfaces of general type&no (Example~\ref{exa:norationalgeneraltype}, Theorem~\ref{thm:genericnorationalcurves}) \\
			& &finitely many? (Example~\ref{exa: fermat surface}, Conjecture~\ref{conj:lang}) \\
			&& \\\hline
		\end{tabular}
		\vspace{1ex}
		\caption{Rational curves on minimal surfaces}
		\label{table: summary on curves over surfaces}
	\end{table}

    \medskip
    
	In Sections~\ref{sec: rat and rul} and~\ref{sect: complex tori}, we see that for several classes of surfaces (namely rational, ruled, tori, hyperelliptic and Kodaira surfaces) the answer to the question of existence of rational curves is complete: either they are covered by them (first two classes) or there are none at all (for the last three). 
	If we move to other classes of surfaces, it is not possible to give a definitive answer to the question. 
	For instance, in Section~\ref{sec: K3}, we give examples of K3 surfaces with either zero, finitely many or infinitely many curves. 
	Moreover, properly elliptic surfaces do also manifest a behavior similar to that of K3 surfaces (see Section~\ref{sec: prop ell}), and it is indeed difficult to determine whether they admit or not rational curves.
	Also for the general type surfaces there are no definitive results. In Section~\ref{sec: gen typ}, we present examples of both existence and non existence of rational curves on them; then we turn our discussion to the Lang Conjecture, that, in a rough form, claims that the number of rational curves in a general type surface must be finite. 
    Finally, in Section~\ref{sect: class VII} we collect the known examples of class VII surfaces and show how the existence of curves distinguishes them.

\subsection*{Notation}
	We denote by $ S $ a connected compact complex non-singular surface and by $ T_S $ its holomorphic tangent bundle.
	We also denote by $ K_S $ the canonical bundle of $S$, namely the line bundle $ \Lambda^2 T_S^* $ of holomorphic $2$-forms on $S$.
	
	 Consider the space $ H^{0}\left(S,K_S^{\otimes m}\right) $ of global pluricanonical forms, where $m\geq1$.
	 If $ {H^{0}\left(S,K_S^{\otimes m}\right)=0} $ for all $m\geq1$, then we say that the \emph{Kodaira dimension} of $S$ is $-\infty$, and we write $\operatorname{Kod}(S) = -\infty$. 
	 If this is not the case, then we let
	\begin{equation*}
		\operatorname{Kod}(S) = \limsup\limits_{m \to \infty} \frac{\log \dim H^{0}\left(S,K_S^{\otimes m}\right)}{\log m}.
	\end{equation*}
	
	The \emph{algebraic dimension} of $ S $, which is denoted by ${a}(S)$, is the transcendence degree over $\C$ of the field of meromorphic functions of $ S $.
 
    Denote by $\chi(S)$ the Euler characteristic of $S$, and by $\chi(S,E)$ the Euler characteristic of a holomorphic vector bundle $E \to S$.
    For $ 0 \le k \le 4 $, the $ k $-th Betti number of $S$ is denoted by $ b_{k}(S) $ or simply $b_k$ if the surface is understood.
	
	Finally, for $ 0 \le p,q \le 2 $, let $ h^{p,q} $ denote the $(p,q)$-th \emph{Hodge number} of $S$, namely the dimension of the $ (p,q) $-Dolbeault cohomology space $ H^{p,q}(S) \cong H^q \left( S,\Lambda^p T_S^* \right) $.
	We arrange the Hodge numbers in the \emph{Hodge diamond} of $S$ as follows
	\begin{equation*}
	    \begin{matrix} 
&  & h^{0,0} & &  \\
& h^{1,0} &  & h^{0,1} &  \\
h^{2,0} &  & h^{1,1} &  & h^{0,2} \\
& h^{2,1} &  & h^{1,2} &  \\
&  & h^{2,2} & &
\end{matrix} 
	\end{equation*}
	The number $h^{0,1}$ is called the \emph{irregularity} of $S$ and is also denoted by $q(S)$.

\section*{Acknowledgements}
The idea to write this survey took place during the conference ``Cohomology of Complex Manifolds and Special Structures - II". The authors are very thankful to the organizers for their warm hospitality and to all the participants for the stimulant environment. 
The authors would also like to thank Daniele Angella and Alexandra Otiman for the help provided during this work and for useful suggestions on the topics of this survey. Many thanks also to the anonymous Referee for the helpful comments and useful suggestions.

\section{Rational and ruled surfaces}\label{sec: rat and rul}
    Understanding the existence of rational curves for class VII$^+_0$ surfaces has proved to be a very challenging problem. On the contrary, for the other surfaces with $ \operatorname{Kod}(S) = -\infty $ this problem is well understood.
    Indeed, rational and ruled surfaces not only admit rational curves, but they are \emph{covered} by them. As we point out in Theorem~\ref{thm:coveringbyrationalcurves}, they turn out to be the \emph{only} surfaces where this is possible. We conclude this section discussing the minimal model program and the \emph{abundance conjecture} (see \cite{KollarMori98} for more details).
    
    \medskip
    
    Algebraic surfaces of Kodaira dimension $-\infty$ were classified at the beginning of the twentieth century by the Italian school of algebraic geometry. They can be divided into three classes:
\begin{itemize}
    \item The projective plane $\P^2$;
    \item The Hirzebruch surfaces $\Sigma_n$, defined as the projective bundles $\P \big(\O_{\P^1}\oplus \O_{\P^1}(-n)\big)$ over $\P^1$;
    \item The projective bundles over a curve of genus $g\geq 1$.
\end{itemize}
In the second item we require $n\neq 1$, otherwise the surface would not be minimal \cite[Chapter~V, Proposition~4.2]{BHPV}.
According to the classification above, we see that these surfaces are covered by rational curves. Furthermore, thanks to the following theorem, we see that they are \emph{the only ones} with this property.
	
\begin{thm}\label{thm:coveringbyrationalcurves}
A {projective} surface $S$ covered by rational curves must have Kodaira dimension~${-\infty}$.
\end{thm}
\begin{proof}
Suppose that there exists a family $\mathcal{C}\subseteq S\times B$ of rational curves parametrized by some variety $B$ such that the morphism $\mathcal{C}\to S$ is dominant, then there exists a dominant rational map
$C\times \P^1\rat S$ where $C$ is a curve, hence $-\infty = \text{Kod}(\P^1\times C) \geq \text{Kod}(S)$. This implies that $\text{Kod}(S)=-\infty$.
	\end{proof}

In spite of the fact that for algebraic surfaces of Kodaira dimension $-\infty$ the existence of rational curves is a straightforward consequence of the classification, it is worth pointing out that it is not clear how to produce rational curves just by knowing that the Kodaira dimension is $-\infty$.
One of the goals of the \emph{minimal model program} is to produce rational curves just by looking at the positivity properties of the canonical divisor. 
The theory in fact was anticipated by the famous \emph{Frankel's conjecture}, which aimed to characterize the projective space as the only complex manifold with positive bisectional curvature. 
Its algebraic formulation asked for positivity properties of the tangent bundle, and it is known as \emph{Hartshorne's conjecture}. Both conjectures were solved independently by Mori \cite{Mori79} and Siu--Yau \cite{SY80}; although completely different, both methods strongly rely on the existence of rational curves. As a matter of fact, the following result is an essential ingredient in Mori's proof.
	
\begin{thm}[\cite{Mori79}]\label{thm:mori}
Let $X$ be a projective algebraic variety and assume that $K_X\cdot C <0$ for a smooth curve $C$ in $X$. Then for all $x\in C$ there exists a rational curve passing through $x$.
\end{thm}

Minimal rational surfaces are the extremal case of the theorem above: their canonical divisor is \emph{anti-ample}, hence intersects every curve negatively, therefore Theorem~\ref{thm:mori} applies proving that they are covered by rational curves. 

Higher dimensional manifolds with anti-ample canonical divisor are called Fano manifolds, and higher dimensional manifolds covered by curves are called \emph{uniruled manifolds}. By Theorem~\ref{thm:mori} a Fano manifold is uniruled, but the converse is not true: a ruled surface has not ample anticanonical divisor. The following version of the so-called \emph{abundance conjecture} aims to characterize uniruledness by the Kodaira dimension.

\begin{conj}[\cite{Kaw85}]
    Let $X$ be a complex projective manifold. Then $X$ is uniruled if and only if $\mathrm{Kod}(X)=-\infty$.
\end{conj}

We have just seen that this conjecture is true in dimension~$2$, while in dimension~$3$ it is a consequence of deep results due to Miyaoka (see \cite{Miy85}). In higher dimensions, the conjecture is still open but there are several partial results, see \cite{BDPP} and references therein.

\section{K3 and Enriques surfaces}\label{sec: K3}
	In this section, we cover the case of K3 and Enriques surfaces, which are classes (5) and (6) of Table~\ref{tab: Enr-Kod classification} respectively. 
    For these classes, the answer to the question of the existence of rational curves is rather subtle. For instance, we give examples of K3 surfaces with either zero, finitely many, or infinitely many curves.
	
	\subsection{K3 surfaces}
	
These are simply connected surfaces with trivial canonical bundle. Every K3 surface is K\"ahler \cite{Siu83}, with Hodge diamond:
\[
\begin{matrix} 
&  & 1 & &  \\
& 0 &  & 0 &  \\
1 &  & 20 &  & 1 \\
& 0 &  & 0 &  \\
&  & 1 & &
\end{matrix} 
\]
The second cohomology group with integral coefficients of any K3 surface endowed with the intersection pairing is a {lattice} of rank $22$, i.e. a free abelian group with an integral quadratic form. We denote this lattice by $\Lambda_{\rm{K3}}$. Moreover, for every K3 surface $S$ the Hodge structure on $H^2(S,\Z)$ induces an \emph{integral} Hodge structure on $\Lambda_{\rm{K3}}$.
	
Another important property of K3 surfaces is the existence of a unique non-degenerate holomorphic $2$-form $\sigma\in H^{2,0}(S)$. The \emph{Torelli Theorem} \cite{PS71} gives a precise characterization of the surface by means of the $2$-form $\sigma$. For our purposes, we only need part of the theorem known as \emph{surjectivity of the period map}.
	
	\begin{thm}[\cite{PS71}]
		For any Hodge structure on the {\rm K3} lattice $\Lambda_{\rm{K3}}$, there exists a {\rm K3} surface $S$ and a Hodge isometry $H^2(S,\Z)\cong \Lambda_{\rm{K3}}$.
	\end{thm}
	
	\begin{rmk}\label{rem:K3withnocurves}
	    The theorem above gives a way to construct K3 surfaces with few rational curves. For example, if the Hodge structure on $\Lambda_{\rm{K3}}$ has no $(1,1)$-classes, then it cannot have rational curves since they are all algebraic. This shows that there are K3 surfaces with no rational curves.
	\end{rmk}
	We provide another example.
	
	\begin{prop}
		There exists a {\rm K3} surface $S$ such that:
		\begin{enumerate}
			\item $\mathrm{Pic}(S)$ is generated by a line bundle $L$ with $c_1(L)^2 = -2$;\label{propertyk3surface pic generated}
			\item either $L$ or $L^{-1}$ is effective with a unique section that is a smooth rational curve;
			\item there are no other irreducible curves in $S$.
		\end{enumerate}
	\end{prop}
	\begin{proof}
		The first property follows by the surjectivity of the period map. Let $S$ be such a surface and $L$ be the generator of $\mathrm{Pic}(S)$. By Riemann--Roch Theorem we have
		\[
		\chi(S,L) = h^0(S,L) -h^1(S,L) + h^2(S,L)  = c_1(L)^2/2+ 2=  1.
		\]
		Serre duality gives $h^2(S,L) = h^0(S,L^{-1})$, hence either $L$ or $L^{-1}$ have sections and in such case the space of sections is 1-dimensional. 
		Without loss of generality, we can assume that $L$ is effective, then it has a unique section $R$ which must be irreducible by the assumption on the Picard group of $S$. By the genus formula, $R$ has arithmetic genus $0$ and therefore is smooth.
		The last property is a consequence of the Lefschetz Theorem on $(1,1)$ classes.
	\end{proof}
	
	A more explicit example of K3 surfaces with rational curves is given by the \emph{Kummer surface} associated with an abelian surface, given in the example below.
	\begin{exa}\label{exa:kummer}
		Let $X$ be a complex torus of dimension $2$ (cf. Section~\ref{sect: complex tori} for the definition) then $X$ admits an involution $\iota$ given by multiplication by $-1$. The quotient $X/\iota$ is a singular surface that has a resolution given by a K3 surface. If one starts with an abelian surface with no curves, our resulting K3 surface will only have $16$ rational curves, given as the resolution of the fixed points of $\iota$. 
	\end{exa}
	
	We now see that whenever the K3 surface is \emph{projective} there exist infinitely many rational curves. The theorem below was proved by Mukai and Mori \cite{MM83}, although they gave credit to Bogomolov and Mumford. We give a sketch of the proof.
	
	\begin{thm}[\cite{MM83}]\label{thm:bogomolovmumford}
		Let $S$ be a projective \rm{K3} surface, then for every ample divisor $H$ on $S$ there exists a rational curve $C$ in $S$ such that $[C] = H$.
	\end{thm}
	\begin{proof}
		The fundamental idea in the proof is to exhibit a K3 surface with two smooth curves $C_1$ and $C_2$ such that $[C_1+C_2] = H$ and $[C_i] \neq H$ for $i=1,2$. Then deform $C_1+C_2$ so that the class becomes irreducible in nearby fibers. The first part of the argument uses the Kummer surface associated with the product of two elliptic curves, with the $C_i$ being exceptional curves of the resolution of the fixed points as in Example~\ref{exa:kummer}.
	\end{proof}
	
	Theorem~\ref{thm:bogomolovmumford} admits several generalizations. It is proved in the same paper \cite{MM83} that for \emph{general} K3 surfaces, there exist \emph{infinitely}-many rational curves. Later it was also proved for other types of surfaces, for example, K3 surfaces having an elliptic fibration \cite{BT00}. The existence of rational curves on any projective K3 surfaces was proved very recently.
	
	\begin{thm}[\cite{CGL19}]\label{thm:infinityrationalcurvesK3}
		There exist infinitely many rational curves on a projective \rm{K3} surface.
	\end{thm}
	
	\begin{rmk}
		Theorem~\ref{thm:infinityrationalcurvesK3} does not contradict Theorem~\ref{thm:coveringbyrationalcurves}. Indeed, the rational curves in a K3 surface cover only a dense subset.
	\end{rmk}
	
	\subsection{Enriques surfaces}
	This class of surfaces was discovered by Enriques at end of the XIX century seeking to answer a question of Castelnuovo about the characterization of rational surfaces. Every Enriques surface is algebraic and has a K3 surface as $2:1$ étale cover. Hence, specifying an Enriques surface is the same as specifying a K3 surface with a fixed-point-free involution. As a consequence, the existence of rational curves on Enriques surface is the same as the one for K3 surfaces. Thus, we can immediately deduce the following result.
	
	\begin{thm}\label{thm:Enriquesinfinitecurves}
		There exist infinitely many rational curves on an Enriques surface.
	\end{thm}
	
	Notice that for the theorem above we do not really need to apply Theorem~\ref{thm:infinityrationalcurvesK3}. Indeed, it is a well-known fact that every Enriques surface has an elliptic fibration, and, therefore, also its K3 cover. Hence, the main theorem of \cite{BT00} applies.

\section{Complex tori, Hyperelliptic and Kodaira surfaces}\label{sect: complex tori}
In this section, we discuss complex tori, hyperelliptic and Kodaira surfaces (i.e., classes (7), (8), and (9) of Table~\ref{tab: Enr-Kod classification} respectively). We gather them together because none of them admit rational curves (see Proposition~\ref{prop:norationalcurvestorus}).
	
	\subsection{Complex tori}
	A complex torus is a manifold of the form $\C^n/\Lambda$, where $\Lambda\subseteq \C^n$ is a lattice of rank $2n$. All complex tori are K\"ahler, have trivial canonical bundle, and can be either projective or not accordingly to their lattice structure. When complex tori are projective they are also called \emph{abelian varieties}. Since the universal cover of any complex torus is $\C^n$, there is a {topological} obstruction to the existence of rational curves.
	
	\begin{prop}\label{prop:norationalcurvestorus}
		Complex tori do not admit rational curves.
	\end{prop}
	\begin{proof}
		Let $X:=\C^n/\Lambda$ and  $f: \P^1\to X$ be any holomorphic map. Since $\C^n$ is simply connected, there exists a unique holomorphic extension $F:\P^1\to \C^n$ of $f$. Since $\P^1$ is compact, by the maximum principle $F$ must be constant, and therefore $f$ is constant.
	\end{proof}
	
	\subsection{Hyperelliptic surfaces}
	A hyperelliptic surface is a surface whose {\em Albanese map} is an elliptic fibration. 
	This class of surfaces shares similitude with Enriques surfaces because they are quotients of a surface with trivial canonical bundle as they are quotients of abelian surfaces.
	More precisely, a hyperelliptic surface is the quotient $(E_1\times E_2)/G$ where $E_1$ and $E_2$ are elliptic curves, and $G$ is a finite group acting freely on $E_1\times E_2$. Since the quotient map is a covering map, the same proof as in Proposition~\ref{prop:norationalcurvestorus} yields the following.
	
	\begin{cor}\label{cor:norationalhyperelliptic}
		Hyperelliptic surfaces do not admit rational curves.
	\end{cor}
	
\subsection{Kodaira surfaces}\label{sect: kodaira surfaces}
	Kodaira surfaces are smooth compact complex surfaces with Kodaira dimension $ 0 $ and odd first Betti number.
	Hence, Kodaira surfaces are not K{\"a}hler.
	They are divided into two families: \emph{primary} and \emph{secondary}, according to whether $ b_1 $ is $ 3 $ or $ 1 $.
	
	\subsubsection*{Primary Kodaira surfaces}
	These are Kodaira surfaces with trivial canonical bundle and first Betti number equal to $3$. All primary Kodaira surfaces are elliptic fiber bundles over elliptic curves. Their Hodge diamond has the form:
	
	\[
	\begin{matrix} 
		&  & 1 & &  \\
		& 1 &  & 2 &  \\
		1 &  & 2 &  & 1 \\
		& 2 &  & 1 &  \\
		&  & 1 & &
	\end{matrix} 
	\]
	We recall that the universal cover of a primary Kodaira surface $ X $ is isomorphic to $ \mathbb{C}^2 $~\cite[Theorem~19]{Kod64}.
	Moreover, the fundamental group $ \pi_1(X) $ can be identified with a group $ \Gamma $ of affine transformations without fixed points of $ \mathbb{C}^2 $ which acts properly discontinuously on $ \mathbb{C}^2 $.
	We can thus identify $ X $ with the quotient $ \mathbb{C}^2 / \Gamma $.
	
	Similarly to the case of complex tori (cf. Proposition~\ref{prop:norationalcurvestorus}), we deduce the following result.
	\begin{prop}\label{prop: primary kodaira surfaces case}
		Primary Kodaira surfaces do not admit rational curves.
	\end{prop}
	
	\subsubsection*{Secondary Kodaira surfaces}
	These surfaces have first Betti number equal to $1$,
	they have torsion canonical bundle, and their Hodge diamond is of the following form:
	
	\[
	\begin{matrix} 
		&  & 1 & &  \\
		& 0 &  & 1 &  \\
		0 &  & 0 &  & 0 \\
		& 1 &  & 0 &  \\
		&  & 1 & &
	\end{matrix} 
	\]
	From~\cite[V.5]{BHPV} we recall that given a secondary Kodaira surface $ Y $, there is a primary Kodaira surface $ X $ such that $ Y $ is a quotient of $ X $ by a cyclic group of finite order.
	It is the group of automorphisms of $ X $, which acts freely on $ X $. 
 
	It follows that a secondary Kodaira surface has a non-ramified finite covering by a primary one.
	In particular, the universal covering of a secondary Kodaira surface is isomorphic to $ \mathbb{C}^2 $.
	
	Consequently, as in the case of the primary Kodaira surfaces, we deduce the following result.
	\begin{prop}\label{prop: secondary kodaira surfaces case}
		Secondary Kodaira surfaces do not admit rational curves.
	\end{prop}
	
	In general, whenever we have a complex manifold whose universal cover does not admit rational curves (for example $\C^n$), the same argument of Proposition~\ref{prop:norationalcurvestorus} yields the non-existence of rational curves.

\section{Properly elliptic surfaces}\label{sec: prop ell}
	
	Here we study properly elliptic surfaces, which represent the class (10) of Table~\ref{tab: Enr-Kod classification}. For this class is already difficult to determine whether they admit or not rational curves and, as in the case of K3 surfaces, there are families for which there are no rational curves.
	
	\medskip
	
	A complex surface $S$ is called \emph{elliptic} if there exists a morphism $\pi:S\to B$ onto a non-singular curve such that the general member is an elliptic curve. If the Kodaira dimension of $S$ is $1$ we say that $S$ is \emph{properly elliptic}. Properly elliptic surfaces can be either K\"ahler or non-K\"ahler.
	
	A properly elliptic surface might have \emph{singular} fibers. However, in \cite{Kod63} Kodaira classified all the possible singular fibers in an elliptic fibration, giving a very precise description of the possible singularieties.
	
	\begin{thm}[\cite{Kod63}]\label{thm:kodairasingularfibers}
		Let $\pi:S\to B$ be an elliptic fibration. Its singular fibers are either a configuration of rational curves or a multiple of a elliptic curve.
	\end{thm}
	
	Hence, whenever there are singular fibers which are not multiple, the surface has at least one rational curve. The study of rational curves therefore falls into two classes, the K\"ahler and non-K\"ahler properly elliptic surfaces.
	
	\subsection{Properly elliptic Non-K\"ahler surfaces}
	If $S$ is non-K\"ahler then $\chi(S) = 0$ and the first Betti number is odd (see~\cite{Buch99,Lamari99}). In particular, if $g$ is the genus of the base curve $B$, then $b_1(S) = 2g+1$, and therefore $q(S)= g+1$.
	
	The following theorem gives a negative answer to the existence of rational curves in this case.
	
	\begin{thm}[\cite{Bri96}]\label{thm:nonkahlerproperlyelliptic}
		Non-K\"ahler properly elliptic surfaces do not contain rational curves.
	\end{thm}
	\begin{proof}
		We give a sketch of the proof. See \cite[Chapter III, Sections 17-18]{BHPV} for the relevant definitions. It is proved in \cite[Chapter~III, Theorem~18.2]{BHPV} that $\deg(R^1\pi_*\O_S)\leq 0$ and it is zero if and only if all non-singular curves are isomorphic and the singular curves are multiples of elliptic curves. 
		The Leray spectral sequence associated to the fibration induces an exact sequence in cohomology
		\[
		\begin{tikzcd}
			0\ar[r] & H^1(B,\O_B)\ar[r] & H^1(S,\O_S)\ar[r] & H^0(B,R^1\pi_*\O_B)\ar[r] & 0
		\end{tikzcd}
		\]
		Since $q(S) = g+1$, we get that $H^0(B,R^1\pi_*\O_B) \cong \C$ and this implies that $\deg(R^1\pi_*\O_B) = 0$.
	\end{proof}
	
	The theorem above also implies that a non-K\"ahler properly elliptic surface is \emph{isotrivial}, i.e. all the smooth fibers are isomorphic.
	
	\subsection{Properly elliptic K\"ahler surfaces}
	The following proposition, which also applies for the non-K\"ahler case, characterizes all the properly elliptic surfaces with no rational curves on their singular fibers.
	
	\begin{prop}\label{prop:chi=0singularfibers}
		Let $\pi:S\to B$ be a properly elliptic surface, then $\chi(S) = 0$ if and only if all the singular fibers are multiple of elliptic curves.
	\end{prop}
	\begin{proof}
	 Denote by $\Sigma\subseteq B$ the critical values of $\pi$, which is a finite collection of points. Thanks to the properties of the topological Euler characteristic, we have
	 \[
	 \chi(S) = \chi_c\left(S\setminus\pi^{-1}(\Sigma)\right) + \chi(\Sigma),
	 \]
	 where $\chi_c$ stands for the Euler characteristic with compact support. Since $S\setminus\pi^{-1}(\Sigma)$ is a smooth elliptic fibration over $B\setminus \Sigma$ we have that 
	 \[
	 \chi_c\left(S\setminus\pi^{-1}(\Sigma)\right) = \chi_c(B\setminus \Sigma) \cdot \chi(S_b) = 0 ,
	 \]
	 where $S_b$ stands for a smooth fiber, which is an elliptic curve. By the Kodaira classification of the singular fibers we have $\chi(S_b) >0$ for every $b\in \Sigma$ if and only if $S_b$ is not a multiple of an elliptic curve (cf. \cite{Kod63}).
	\end{proof}
	
	\begin{cor}\label{cor:properlyelliptic}
	Let $\pi:S\to B$ be a properly elliptic K\"ahler surface. If $\chi(S)\neq 0$, then $S$ contains at least one rational curve. Moreover, if $g(B)\geq 1$, then $S$ does not contain other rational curves beside the components of the singular fibers of $\pi$.
	\end{cor}
	\begin{proof}
	    By Theorem~\ref{thm:kodairasingularfibers} and Proposition~\ref{prop:chi=0singularfibers} if $\chi(S) \neq 0$ there exists a singular fiber which contains a rational curve. Let $C$ be a rational curve on $S$ that is not contained in a fiber of $\pi$, then the restriction of $\pi$ to $C$ gives a surjective map into $B$, which implies that $g(B) = 0$.
	\end{proof}
	
	\begin{rmk}\label{rmk: elliptic algebraic no curves}
	Corollary~\ref{cor:properlyelliptic} also implies that if $\pi:S\to B$ is a properly elliptic surface with $\chi(S) = 0$ and $g(B)\geq 1$, then it does not have rational curves.
	\end{rmk}
	
	\begin{exa}\label{ex: elliptic algebraic no curves}
	    Corollary~\ref{cor:properlyelliptic} applies to give a very simple properly elliptic surface without rational curves. 
	    We just need to consider the product $B\times E \to B$, where $E$ is an elliptic curve and $g(B)\geq 2$. 
	    Notice that $\chi(B\times E) = 0$.
	\end{exa}
	
	An elliptic surface $\pi:S\to B$ can be viewed as an \emph{elliptic curve} over the field $\C(B)$ of meromorphic functions on $B$. This point of view helps to bring methods in the theory of elliptic curves to study the existence of rational curves on $S$. In particular, when $\pi$ admits a section, the rational curves that map surjectively to $\P^1$ form a group, called the \emph{Mordell--Weil group} of $\pi$. An arithmetic approach via the Mordell--Weil group yields the following result.
	
	\begin{thm}[\cite{Ulm17}]\label{thm: ulmer theorem}
	   Let $\pi:S\to \P^1$ be a very general\footnote{Here very general means outside of a countable number of proper closed subsets in the moduli space.} properly elliptic surface, then $S$ has no rational curves other than the zero section and the components of the singular fibers.
	\end{thm}

\begin{rmk}\label{rmk: infinite curves on elliptic}
    Whenever there exists a rational multisection that is non-torsion on the Mordell--Weil group, then the surface contains infinitely many rational curves, cf. \cite[Section 6]{Ulm17} for details.
\end{rmk}

\section{Surfaces of general type}\label{sec: gen typ}
	
	In this section, we turn to general type surfaces (i.e., class (11) of Table~\ref{tab: Enr-Kod classification}). After giving some theorems and examples concerning the existence or non-existence of rational curves, we focus on the Lang Conjecture \cite{Lang71}, that in a rough form claims that the number of rational curves in a general type surface must be finite.

    \medskip
	
	Surfaces of general type are complex surfaces of Kodaira dimension~2. We can think of them as the 2-dimensional analog of curves of genus greater or equal than 2.  As the name suggests, general type surfaces are in several ways the \emph{most generic} algebraic surfaces. The so-called \emph{geography} of surfaces aims to provide examples of general type surfaces with prescribed Chern numbers; this is an active research area \cite{MendesPardini12,Persson85,RoulleauUrzua15}. 
	
	Every general type surface is projective and its canonical divisor is big and nef. That being said, it is not possible to use Theorem~\ref{thm:mori} to produce rational curves, nevertheless, when the divisor is not ample we can guarantee the existence of at least a rational curve.
	
	\begin{prop}\label{prop: general type K=0}\cite[Exercise 8, page 219]{Debarre2001}
		Let $S$ be a surface of general type and $C$ an irreducible curve on $S$. Then $K_S\cdot C = 0$ if and only if $C$ is a smooth rational curve.
	\end{prop}
	
	From now on, thanks to Proposition~\ref{prop: general type K=0}, we focus on the case of ample canonical divisor. There are very simple examples where it can be shown the non-existence of rational curves, we provide one below.
	
	\begin{exa}\label{exa:norationalgeneraltype}
		Let $C_1,C_2$ be curves of genus greater or equal than $2$, then their product $C_1\times C_2$ is a surface of general type that does not contain rational curves. Indeed, suppose there exists a rational curve $R$ in $C_1\times C_2$ and denote by $\pi_i$ the projection to the $i$-th factor. The curve $R$ cannot be contained in a fiber of $\pi_i$, hence, being the map proper, it must map surjectively to $C_i$, which is a contradiction by Riemann--Hurwitz Theorem.
	\end{exa}
	
	A series of examples of general type surfaces are given by hypersurfaces of degree $d\geq 5$ in $\P^3$. Even in this case, the existence of rational curves is not an easy question. To give a sample of the type of results in this direction we state a theorem due to Clemens.
	
	\begin{thm}[\cite{Cle86}]\label{thm:genericnorationalcurves}
		Let $S\subseteq \P^3$ be a generic smooth surface of degree $d\geq 5$. Then there are no rational curves on $S$.
	\end{thm}
	
	On the other hand, for \emph{specific} general type surfaces, the statement in Theorem~\ref{thm:genericnorationalcurves} no longer holds, as we see in the following example.
	
	\begin{exa}\label{exa: fermat surface}
		Consider the Fermat surface
		\[
		S_d:= \{ x_0^d +x_1^d + x_2^d + x_3^d=0 \} \subset \P^3.
		\]
		For $d\geq 5$ this is a surface of general type. Notice that it is easy to exhibit a family of $3d^2$ lines on $S_d$ and, in fact, Fermat surfaces contain \emph{exactly} $3d^2$ lines. In particular, $S_d$ contains at least $3d^2$ rational curves. However, it is not known whether this is the precise number of rational curves, yet.
	\end{exa}
	
	Fermat surfaces give examples of general type surfaces containing any number of rational curves, so a natural question arises: is it possible to have an \emph{infinite} number of rational curves? This question is thigh-related with the \emph{hyperbolicity} of our surfaces and the so-called Lang conjecture, which we state here in a simplified form.
	
	\begin{conj}[\cite{Lang74}]\label{conj:lang}
		Let $S$ be a surface of general type. There exists at most a finite number of rational curves on $S$.
	\end{conj}
	
	Let us note that the conjecture above needs the general type assumption. Indeed, as we saw before, it is false for rational, ruled, K3, and Enriques surfaces. 
	
	There are two partial results that provide evidence to Conjecture~\ref{conj:lang}.
	
	\begin{thm}[\cite{Lu95}] 
		Let $S$ be a surface of general type. There exists at most a finite number of \emph{smooth} rational curves on $S$.
	\end{thm}
	
	Therefore, this theorem leaves open the problem of counting {singular} rational curves.
	
	One of the most important results related to Conjecture~\ref{conj:lang} is the following theorem due to Bogomolov.
	
	\begin{thm}[\cite{Bog78}]
		Let $S$ be a general type surface with $c_1^2-c_2>0$, then $S$ contains at most a finite number of rational curves.
	\end{thm}
	
	Summarizing, we have seen that a sufficient condition for the existence of rational curves on general type surfaces seems to be unclear (even conjecturally) up to now. 
	
\section{Class VII surfaces}\label{sect: class VII}
    As we saw in the introduction, the classification of class VII surfaces remains the main open problem in the Kodaira classification of surfaces and it has reduced to the \emph{Global Spherical Shell Conjecture} (in short, \emph{GSS Conjecture}).
    
	\begin{defi}[\cite{Kat78}]
		A \emph{global spherical shell} in a surface $S$ is a neighborhood $V$ of $\S^3 \subset \C^2 \setminus \{0\}$ which is holomorphically embedded in the surface $S$ so that $S \setminus V$ is connected.
	\end{defi}
	The global spherical shell conjecture claims that all class VII$_{0}$ surfaces with positive second Betti number have a global spherical shell.
    As a matter of fact, it is known that a class VII surface with positive second Betti number $b_2 > 0$ has at most $b_2$ rational curves and has exactly this number if it has a global spherical shell \cite{Dlo84}.
	Conversely, Dloussky, Oeljeklaus, and Toma \cite{DOT} showed that if a class VII$_{0}$ surface with positive second Betti number $b_2>0$ has exactly $b_2$ rational curves, then it has a global spherical shell. In a series of celebrated papers \cite{Tel05, Tel06, Tel09} Teleman developed a strategy using methods of gauge theory to find rational curves. For this, we need to introduce the definition of a cycle of rational curves.
	
\begin{figure}[h!]
        \centering
        \includegraphics[width=6.2cm]{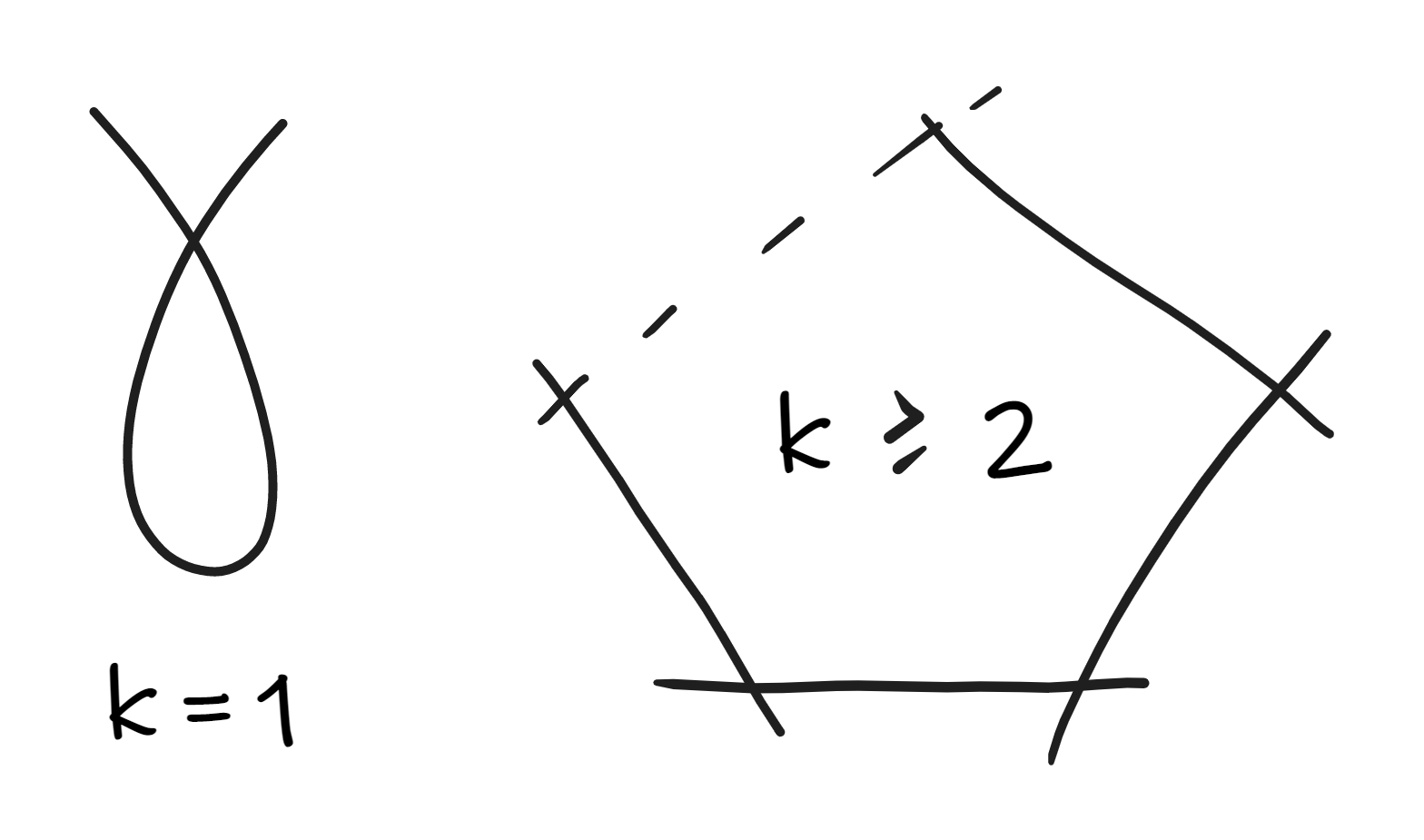}
        \caption{Cycles of rational curves.}
        \label{fig cycle}
    \end{figure}

	\begin{defi}\label{def: cycle}
        A \emph{cycle} of rational curves on a complex surface $S$ is an effective divisor $D\subset S$ which is either a rational curve with a simple singularity or a sum of $k\geq2$ smooth rational curves intersecting as in Figure \ref{fig cycle}.
    \end{defi}
	
	In \cite{Tel05} Teleman proved the existence of a rational curve on class VII$_0$ surfaces with $b_2=1$, verifying the GSS Conjecture for these surfaces. 
    Subsequently, in \cite{Tel06}, Teleman extended his result when $b_2=2$ showing the presence of a \emph{cycle} of rational curves; since the maximum possible of rational curves is two, there are two possibilities for such a cycle: either consists on two smooth rational curves or just one rational curve with a simple singularity. 
    Therefore, for the latter case, one needs to find another rational curve in order to conclude that it has a global spherical shell. 
    There have been attempts to understand whether there actually is a second rational curve.
    For example, in \cite{Brunella11,Brunella11corig} M. Brunella studied class VII$_0$ surfaces with $b_2=2$ and such that they admit a \emph{holomorphic foliation}.
    
    Finally, we remark that for class VII$_0$ surfaces with $b_2 > 2$ the existence of rational curves is still widely open, although progress has been done.
    See for example the contributions by G. Dloussky \cite{Dloussky06, Dloussky21}.
    It is worth mentioning that there are further works trying to understand whether a surface of class VII$_0$ has rational curves or not by adding more assumptions.
    For instance, in \cite{Brunella14} M. Brunella studied this problem with the extra condition of the existence of an automorphic Green function on a cyclic covering.
    We shall also mention the works of M. Pontecorvo and A. Fujiki \cite{Pontecorvo17, FP19}, and V. Apostolov \cite{Apostolov01} which seek to understand class VII surfaces through the study of bi-Hermitian metrics.\\

    We highlight that a proof of the GSS conjecture would lead to a complete classification of the type VII surfaces. 
    Indeed, the surfaces with a global spherical shell are reasonably well understood since they are \emph{Kato surfaces} \cite{Kat78}. 
    As a matter of fact, a construction method and thus a large class of examples have been introduced by Kato \cite{Kat78}, and at present, these are the only known minimal surfaces of class VII with $b_2 > 0$.

\subsection{Kato surfaces}\label{subsection:katosurfaces}
    The sphere $\S^1$ in a torus $\S^1\times\S^1$ is a trivial example of a sphere in a manifold that does not disconnect it. Similarly, the following construction is the most natural way to obtain a manifold with a global spherical shell. Define the ball and the sphere of radius $r>0$ in $\C^2$ as
    $$ B_r := \{ z\in \C^2 \;\text{s.t.}\; |z|<r\} ,\;\; S_r:=\p B_r . $$
    Consider $\overline{ B_1} \setminus B_{1/2}$, and identify the two boundaries of this space, which are $S_1$ and $S_{1/2}$. The image of a small neighborhood of $S_1$ in the manifold so obtained gives us a global spherical shell. 

    \begin{rmk}
        The manifold constructed above is the Hopf surface $\S^3\times\S^1$, see Section~\ref{ssec: class by curves}.
    \end{rmk}

\begin{figure}[h!]
        \centering
        \includegraphics[width=8cm]{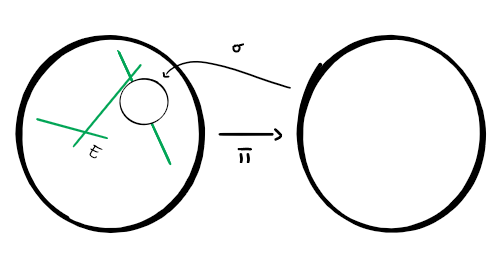}
        \caption{Construction of Kato surfaces.}
        \label{fig Kato1}
    \end{figure}

    To construct a Kato surface we just need to add an intermediate step to this construction. Namely, we consider a modification of $B_1$ at finitely many points, $\pi: B \longrightarrow B_1$. We then remove a small open ball from $B$, $\sigma: B_1 \hookrightarrow B$. Finally, the compact complex surface $X$ is obtained by identifying the two boundaries of $\overline{B}\setminus\sigma(B_1)$ (see Figure \ref{fig Kato1}).

    Notice that we obtain a minimal surface if and only if the modification $\pi: B\longrightarrow B_1$ is a sequence of blow-ups $\pi = \pi_n \circ \cdots \circ \pi_1$ such that $\pi_j$ is the blow-up of a point $p_{j-1}$ in the exceptional divisor of $\pi_{j-1}$ and $\pi_1$ is a blow-up in a point $p_0\in B_1$ such that $\sigma(p_0)$ is in the exceptional locus of $\pi$.
    This implies that $p_0$ is a fixed point for $\pi\circ\sigma$. 
    In other words, to obtain a minimal surface out of this construction we perform successive blow-ups of the unit ball $B_1$ starting from the origin and continuing to blow up a point on the last constructed exceptional curve until, in order to get a compact complex surface, we take the quotient by the biholomorphism $\sigma$.
    As a matter of fact, this construction has a trivial case in which it reduces to a modification of the Hopf surface: this happens when $\sigma(p_0)$ does not meet the exceptional locus.
    However, in any case, Kato surfaces are deformations of modifications of Hopf surfaces (\cite{Nak90}, see also the next section).
    
    \medskip
    
    The Kato surfaces are the only known examples of class VII$_0^+$ surfaces and it was conjectured by Nakamura \cite[Conjecture~5.5]{Nak84b} that they represent the whole class. 
    Indeed, they produce examples of class VII$_0$ surfaces for any $b_2\geq 0$. 
    The second Betti number of a minimal Kato surface is exactly the number of blow-ups of the modification $\pi$ and it agrees with the number of rational curves on it. In particular, the above process constructs $n$ rational curves $D_1, \ldots , D_n$ of self-intersection less or equal to $-2$ which are the only rational curves in the surface. 
    
    \medskip
    
    We now look at the simple example of $b_2=1$.
    We use it to understand how the exceptional divisor becomes a rational curve with self-intersection different from $-1$ after the identification. Obviously, we shall suppose that $\sigma(0)\in E$, where $E$ is the exceptional divisor of the blow-up in $0$ of $B_1$ given by $\pi$. 
    
\begin{figure}[h]
        \centering

            \includegraphics[width=12cm]{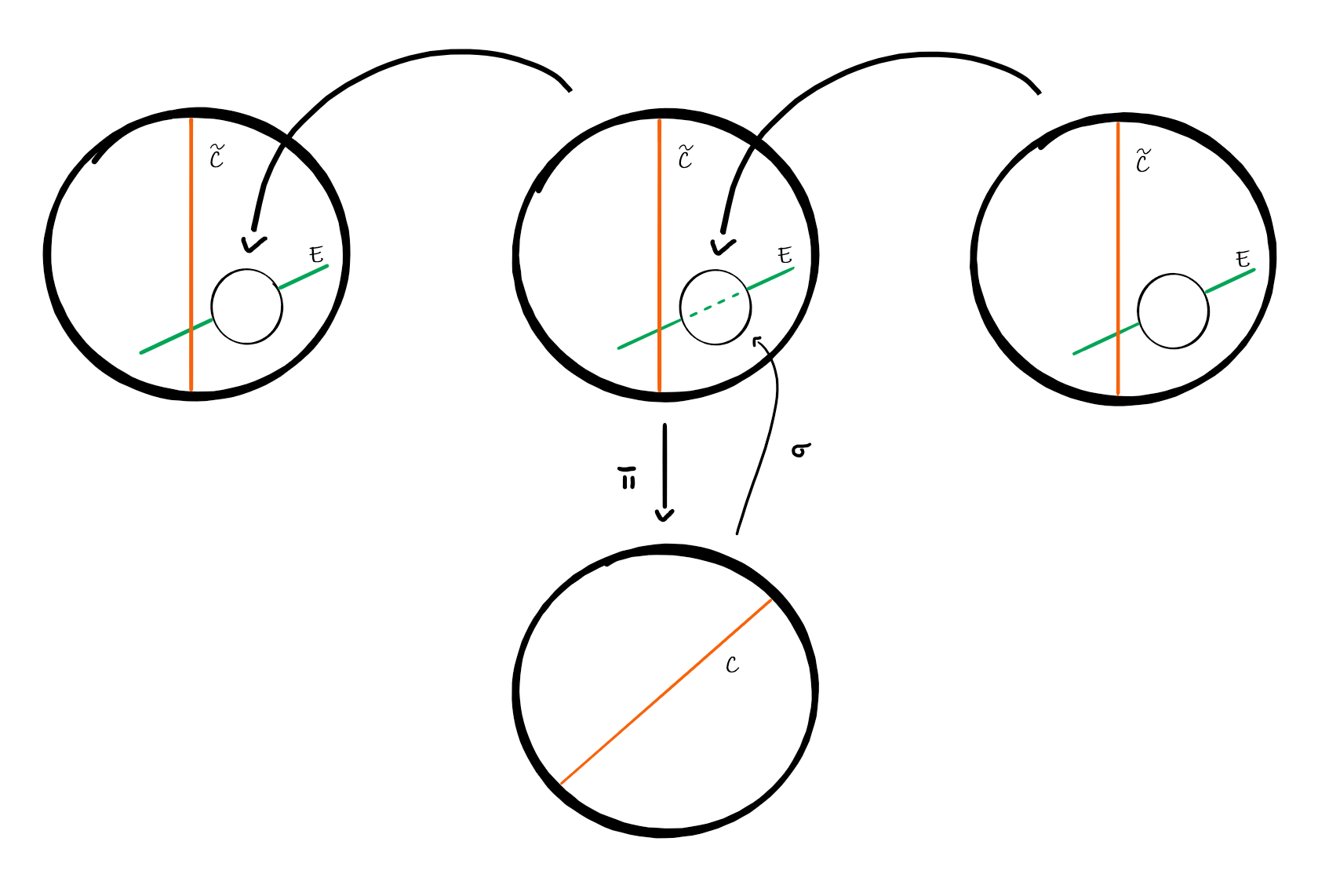}
        \caption{Universal cover of a Kato surface.}
        \label{fig: Kato2}
    \end{figure}

    We define $\mathcal{C}=\sigma^{-1}(E)$ and $\widetilde{\mathcal{C}} = \bar{\pi^{-1}(\mathcal{C}\setminus\{0\})}$. Since the identification is made by $\sigma\circ\pi$, the universal cover of this Kato surface is made by infinitely many copies of $\bar B \setminus \sigma(B_1)$ glued so that $\widetilde{\mathcal{C}}$ of one copy is completing $E$ of the next one (see Figure \ref{fig: Kato2}). 
    
\begin{figure}[h!]
        \includegraphics[width=6cm]{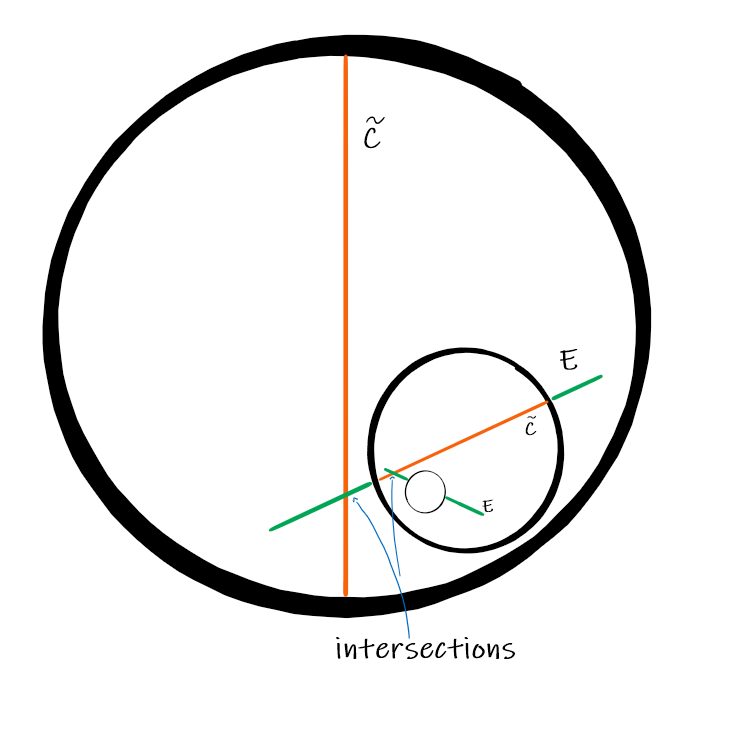}
        \caption{Two copies of $\bar B \setminus \sigma(B_1)$ glued through $\sigma\circ\pi$.}
        \label{fig Kato3}
    \end{figure}

    In this simple example, we are also supposing that $\sigma(0) \notin \pi^{-1}(\mathcal{C})$. In this way, the intersection between these two lines is $ \widetilde{\mathcal{C}}\cdot E = 1$.
    Thus the universal cover contains an infinite chain of $(-2)$-curves that transversely intersect one to each other at one point.
    We now show that they have self-intersection $-2$.
    This is because they have the same self-intersection of $\widetilde{\mathcal{C}}$ which actually is $-2$.
    We perform our computations in the union of two copies of $\bar B \setminus \sigma(B_1)$ (as in Figure \ref{fig Kato3}) so that the ambient space is compact.
    
    First of all, $E^2=-1$ implies $\mathcal{C}^2=-1$ since $\mathcal{C}=\sigma^{-1}(E)$.
    Then $(\widetilde{\mathcal{C}} + E)^2=-1$ since $\widetilde{\mathcal{C}} + E = \pi^{-1}(\mathcal{C})$.
    Hence,
    $$\widetilde{\mathcal{C}}\cdot\widetilde{\mathcal{C}} = (\widetilde{\mathcal{C}} + E)^2 - E^2 -2E\cdot\widetilde{\mathcal{C}} = -1 +1 -2 = -2.$$

    Finally, we consider a double cover of our Kato surface (Figure \ref{fig Kato3} with quotients of the boundaries). Here the pull-back of our rational curve is given by two $(-2)$-curves which intersect twice. 
    
    We call $\widetilde{E}$ the rational curve on the Kato surface and define the divisors $E_1$ and $E_2$ on the double cover as $\pi^{-1}(\widetilde{E})=E_1 + E_2$, where $\pi$ is the projection of the double cover. We have,
    $$ \widetilde{E}\cdot\widetilde{E} = (E_1 + E_2)^2 = E_1^2 + 2\,E_1\cdot E_2 + E_2^2 = -2 + 2\cdot2 -2 = 0 \,$$
    hence, after the identification, the exceptional divisor becomes a rational curve of zero self-intersection.

\subsection{Classification by curves}\label{ssec: class by curves}
    In this section, we give an overview of the classification of the known families of surfaces in class VII. 
    In particular, we look at {\em Enoki surfaces, parabolic, hyperbolic} and {\em half Inoue surfaces}, {\em Inoue--Bombieri surfaces}, and {\em Hopf surfaces}. 
    In the spirit of this survey, we recall the results of Inoue, Kato, and Nakamura to see how these families of surfaces are characterized by the existence of curves.\\
    
    Hopf manifolds are compact complex non-K{\"a}hler manifolds obtained as quotients \mbox{$\left(\C^n \setminus \{0\}\right) \slash \Z$}, where $\Z$ is acting freely by holomorphic contractions.
    Since $\Z$ is acting by contraction, these manifolds are all diffeomorphic to $\S^{2n-1}\times\S$.
    In the case of complex dimension two, these are called {\em primary Hopf surfaces}, to distinguish them from the {\em secondary Hopf surfaces}.
    The latter are quotients $\left(\C^2 \setminus \{0\}\right) \slash G$ where $G$ acts properly discontinuously on $\C^2 \setminus \{0\}$ but is not infinite cyclic.

    There is a normal form for the contraction \cite[Theorem 1]{MR196769}, which in appropriate coordinates can be written as
    $$ \gamma(x,y) =  (ax + \lambda y^n, by), $$
    where $\gamma$ is the generator of $\Z$, $a,b \in \C$ are such that $0<|a|<1 $, $0<|b|<1$ and either $\lambda =0$ or $a = b^n$. The image of the $x$-axis gives an elliptic curve on these surfaces. Moreover, whenever $\lambda$ vanishes the image of the $y$-axis gives a second elliptic curve on them. We have the following theorem due to Kato.
    
    \begin{thm}[Kato, see (5.2) of \cite{Nak84}]\label{thm: hopf surf}
        Let $S$ be a minimal surface of class \rm{VII} with no meromorphic functions and exactly two elliptic curves then it is isomorphic to a Hopf surface.
    \end{thm}
    
    In \cite{Ino74} Inoue introduced some surfaces in class VII, known as \emph{Inoue--Bombieri} surfaces, with $b_1=1$, $b_2=0$ and no curves by taking quotients of $\H\times\C$, where $\H$ is the upper complex half-plane.
    In \cite{Ino77}, taking different quotients of $\H\times\C$ Inoue obtained other class VII surfaces with $b_1=1$ and $b_2>0$ (with curves!).
    These are known in the literature as {\em parabolic, hyperbolic} and {\em half Inoue surfaces}, see also \cite[(1.1),~(1.3), and~(1.6)]{Nak84} for the precise definitions. 
    We have the following result of Inoue about their properties.    
    \begin{thm}[\cite{Ino77}] The followings hold.
        \begin{itemize}
            \item Parabolic Inoue surfaces have a cycle $C$ of rational curves such that $C^2=0$ and an elliptic curve;
            \item Hyperbolic Inoue surfaces have two cycles $C_1$ and $C_2$ of rational curves and $b_2 = n+m$ where $n$ and $m$ are the numbers of (possibly singular) rational curves of $C_1$ and $C_2$ respectively;
            \item Half Inoue surfaces have a unique cycle $C$ of rational curves with $C^2=-b_2$.
        \end{itemize}
    \end{thm}
    Moreover, the following result by Nakamura shows how the existence of rational curves characterizes these three families.
    \begin{thm}[\cite{Nak84}]
        Given $S$ a minimal surface in class \rm{VII},
        \begin{itemize}
            \item if $S$ has a cycle of rational curves and an elliptic curve then it is isomorphic to a parabolic Inoue surface;
            \item if $S$ has two cycles of rational curves then it is isomorphic to a hyperbolic Inoue surface;
            \item if $S$ has a cycle of rational curves, $C$, such that $b_2=-C^2$ then it is isomorphic to a half Inoue surface.
        \end{itemize}
    \end{thm}
    
    We recall that parabolic Inoue surfaces are a particular case of {\em Enoki surfaces}.
    These were introduced in \cite{Eno81}, and are precisely those class VII$_0$ surfaces with a cycle of rational curves $C$ such that $C^2=0$.
    An Enoki surface admitting a compact curve $E$ which is not rational smooth is a parabolic Inoue surface; as a matter of fact, in this case $E$ must be smooth elliptic with $E^2 = -b_2$.\\

    Kato \cite{Kat79} and Dloussky \cite{Dloussky88,Dlo-Koh} proved that both the Enoki and Inoue surfaces have a global spherical shell, hence they are all Kato surfaces (with $n=b_2$ rational curves). Thus, we can highlight how they arise from the Kato construction. 
    In the language of the previous subsection, the hyperbolic and half Inoue surfaces occur when each blown-up point is at the intersection with a previously created exceptional curve; as we saw, these surfaces are hyperbolic Inoue when they have two cycles and the only other case is that of their $\mathbb{Z}_2$-quotients which are half Inoue surfaces. On the other hand, when every blow-up occurs at a generic point of the previously created exceptional curve, we get $D^2_i = -2$ for each $i = 1,\ldots,n$ and the divisor $D_1 + \cdots + D_n$ is a cycle $C$ of rational curves of self-intersection $C^2 = 0$, hence, in this case, the surface is an Enoki surface. 
    
    Enoki, half, and hyperbolic surfaces are known as the {\em extreme} cases in the terminology of Dloussky \cite{Dlo84}. 
    He introduced an index invariant\footnote{In literature the Dloussky index is usually denoted by $\sigma$. 
    Since we already used $\sigma$ in the Kato surfaces construction, we replaced it with $i_{Dl}$.} given by
    $$ i_{Dl}(S) := - \sum_{j=1}^n D_j^2 \;, $$
    which is bounded between $2b_2$ and $3b_2$.
    Enoki surfaces correspond to $i_{Dl}(S)=2b_2$, while half and hyperbolic surfaces correspond to $i_{Dl}(S)=3b_2$.  
    All other surfaces are called {\em intermediate Kato surfaces} since they have $ 2b_2 < i_{Dl}(S)< 3b_2$. These are obtained by a mixed procedure in which there are generic blow-ups as well as blow-ups at intersection points.
    They have only one cycle composed of $b_2 - N$ rational curves and $N \, (\neq0)$ branched curves (see for example the following Figure~\ref{fig: kato branch}).
    
    \begin{figure}[h!]
        \includegraphics[width=5.2cm]{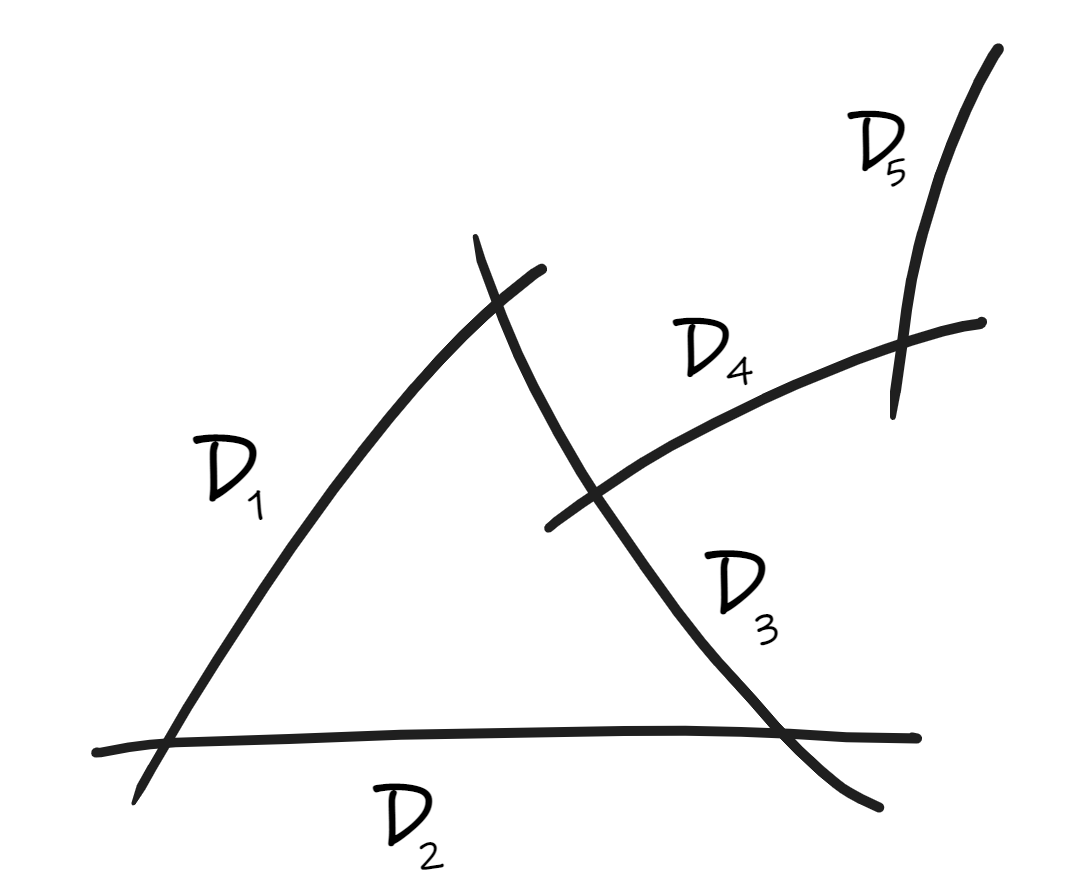}
        \caption{A cycle of three rational curves with two branches}
        \label{fig: kato branch}
    \end{figure}
    
    We can finally summarize all the discussion of this Section~\ref{ssec: class by curves} in the following Table~\ref{table: class VII}.
    In particular, here we collect all the know families of class VII surfaces, and we highlight how they can be distinguished by the existence of curves. 
    
    \begin{table}[h!]
        \centering
        \begin{tabular}{ c|c } 
             curves & surface \\ [0.5ex] 
            \hline \hline 
               no curves & Inoue--Bombieri surface \\ [0.5ex]
                 2 elliptic curves     &  Hopf surface     \\ [0.5ex] 
                2 cycles of rational curves & hyperbolic Inoue surface \\[0.5ex] 
                a cycle $C$ with $C^2=-b_2$ & half Inoue surface  \\[0.5ex]
                a cycle and $N\neq0$ branched curves & intermediate Kato\\[0.5ex]
                      a cycle $C$ with $C^2=0$    & Enoki surface    \\[0.5ex]
                a cycle and an elliptic curve & parabolic Inoue surface \\ [0.5ex] 
                no cycles and an elliptic curve &  Hopf surface   \\[1ex]
        \end{tabular} 
        \caption{List of the known families of minimal class VII surfaces by means of their curves.}
        \label{table: class VII}
    \end{table}
    
    Now, given $S$ a class VII$_{0}$ surface with a cycle of rational curves $C$, it can be proved that $H^2(S,T_S (-\log C))=0$, see \cite[Theorem 1.3]{Nak90}. 
    This means that there are no obstructions to deform $S$ by keeping $C$ of Hodge type (cf. \cite[Theorem 1.4]{Nak90}).
    Hence we can construct a deformation of $S$ that takes $C$ to a non-singular curve, which must be elliptic since $C^2 = 0$. Thus we obtain a deformation whose general fiber is a surface with an elliptic curve.
    Looking at Table~\ref{table: class VII} this must be isomorphic to a blown-up parabolic Inoue surface or generically a blown-up primary Hopf surface, see \cite[Theorem 12.3]{Nak84}. 
    
    Since, for an arbitrary class VII$_0$ surface, having $b_2>0$ rational curves implies the existence of a cycle of rational curves, the above argument enabled Nakamura to prove the following result.
    
    \begin{thm}[\cite{Nak90}]\label{thm: naka90}
        Given a class \rm{VII$_{0}$} surface $S$ with exactly $b_2$ rational curves ($b_2>0$), we can construct a complex analytic family $\pi: \mathcal{S} \rightarrow B_1$ of small deformations of $S$ such that 
        \begin{itemize}
            \item $S=\pi^{-1}(0)$;
            \item $\mathcal{S}_t = \pi^{-1}(t), \; t\neq0$, is a modification of a primary Hopf surface. 
        \end{itemize}
    \end{thm}
    
\begin{figure}[h]
        \includegraphics[width=9cm]{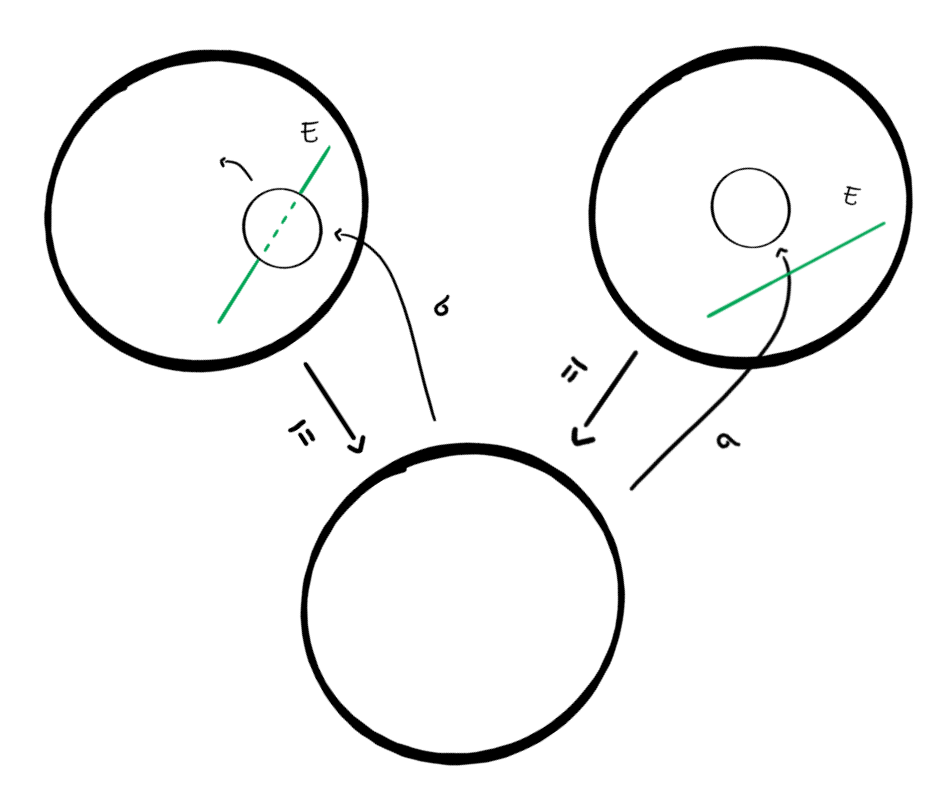}
        \caption{Deformation of Kato surfaces.}
        \label{fig: deform Kato}
    \end{figure}

    We can give an interpretation of this theorem by looking at the Kato surfaces. Indeed, the surfaces in the statement have a global spherical shell, since they have exactly $b_2$ rational curves. We have also seen that modifications of Hopf surfaces are obtained as a trivial case of Kato surfaces precisely when $\sigma(B_1)$ does not meet the exceptional locus. Moreover, in any case, we can perform a translation of $\sigma(B_1)$ in $B$ in order to reduce to the trivial case (Figure \ref{fig: deform Kato}).
    
    We remark that Kato's construction may be performed in all dimensions giving Kato manifolds, see \cite{IOP} for details. Moreover, Kato proved that all compact complex manifolds of complex dimension $n\geq 2$ containing a global spherical shell are constructed in this way. The following result of Kato (which generalizes Theorem~\ref{thm: naka90}) holds.
    
    \begin{thm}[\cite{Kat78}]
        Suppose that a compact complex manifold $X$ of complex dimension $n\geq 2$ contains a global spherical shell. Then we can construct a complex analytic family $ \pi:\mathcal{X}\rightarrow B_1 $ of small deformation of $X$ such that
        \begin{itemize}
            \item $X = \pi^{-1}(0)$;
            \item $X_t = \pi^{-1}(t), \; t\neq 0$, is biholomorphic to a compact complex manifold which is a modification of a Hopf manifold at finitely many points. 
        \end{itemize}
    \end{thm}

\subsection{Analytic approach to the GSS Conjecture}
    We want to briefly recall the link between the existence of curves (with negative self-intersection) on surfaces and the evolution of geometric flows of metrics. We shall outline that some of these flows, such as the pluriclosed flow (see below), have been introduced as analytic tools with the aim of detecting these divisors.
	
	Eyssidieux, Guedj, Song, Tian, Weinkove and Zeriahi studied the analytic minimal model program \cite{EGZ, SoT12, SoT17, SW13, SW14}, seeking to understand the algebraic minimal model program through the singularities of K{\"a}hler--Ricci flow. It was shown in \cite{SW13} and \cite{SoW13} that the K{\"a}hler--Ricci flow on a K{\"a}hler surface contracts $(-1)$-curves in the sense of Gromov-Hausdorff and converges smoothly outside of the curves. Indeed, thanks to the Nakai--Moishezon criterion of \cite{Buc00} and \cite{Lam99} it can be proved that, if the maximal existence time of the flow is finite, then either the volume of the surface goes to zero, or the volume of a curve of negative self-intersection goes to zero. 
	In \cite{TW}, the same behavior was proved to occur also for the Ricci flow on non-K{\"a}hler surfaces with the extra assumption of starting from a pluriclosed metric (which means that $dd^c\omega=0$ where $\omega$ is the K{\"a}hler form associated to the metric). 
	
	We shall remark that Ricci flow may not be the most natural tool for studying geometric evolution on non-K{\"a}hler surfaces (or manifolds in general). Indeed, for non-K{\"a}hler metrics the Ricci flow will not even preserve the condition that a metric is Hermitian. With this idea in mind in \cite{ST10} the authors introduced a parabolic flow of Hermitian metrics which preserves the pluriclosed condition (hence it was named {\em pluriclosed flow}). 
	They also conjectured that it should have the same regularity of the K{\"a}hler--Ricci flow, meaning that the flow should exist until either the volume collapses or it becomes singular on an effective divisor with negative self-intersection, see Section~5 of \cite{ST13} for details. 
	If this conjecture holds, then any class VII$_0^+$ surface should contain an irreducible effective divisor of non-positive self-intersection \cite[Theorem 7.1]{ST13}. 
	Thus, the proof of this conjecture would lead to an analytic proof of the global spherical shell conjecture for $b_2=1$. 
	As a matter of fact, by general theory (\cite[Lemma 2.2]{Nak84}) there would be only two possible cases: the curve is either a rational curve or an elliptic curve.
	If the curve is elliptic, the manifold is known by \cite{Nak84,Eno81}. 
	In the other case, it would contain a global spherical shell.
	
	We refer to \cite{Streets} for further analysis of the conjectural behavior of the pluriclosed flow on class VII surfaces, which seems to be a promising direction for the general case $b_2>0$.

\end{document}